\documentclass{amsart} 
\usepackage{latexsym,amssymb,pstricks}

\newtheorem{thm}{Theorem}[section]
\newtheorem{lem}[thm]{Lemma}
\newtheorem{prop}[thm]{Proposition}
\newtheorem{cor}[thm]{Corollary} 
\newtheorem{de}[thm]{Definition} 
\newtheorem{rem}[thm]{Remark}

\newcommand{\BZ}{{\mathbb{Z}}}

\newcommand{\BO}{{\mathcal{O}}}
\newcommand{\BD}{{\mathcal{D}}}
\newcommand{\BV}{{\mathcal{V}}}
\newcommand{\BA}{{\mathcal{A}}}

\newcommand{\BS}{{\mathcal{S}_p}}
\newcommand{\BSplus}{{\mathcal{S}^+_p}}

\newcommand{\BSf}{{\mathcal{S}_5}}
\newcommand{\BSfplus}{{\mathcal{S}^+_5}}

\newcommand{\Vp}{{V_p}}
\newcommand{\Zp}{{Z_p}}

\newcommand{\Si}{{\Sigma}}

\newcommand{\Gtwo}{
\psset{unit=.5cm}
\begin{pspicture}[.4](-.5,0)(1.5,1)
\pscircle(0,.5){.35}
\pscircle(1,.5){.35}
\psline(.35,.5)(.65,.5)
\end{pspicture}
}
\newcommand{\Gn}{
\psset{unit=.5cm}
\begin{pspicture}[.4](-.5,-.2)(5,1)
\pscircle(0,.5){.35}
\pscircle(1,.5){.35}
\psline(.35,.5)(.65,.5)
\psline(1.35,.5)(1.65,.5)
\psline(3.35,.5)(3.65,.5)
\pscircle(4,.5){.35}
\put(2,.45){$\cdots$}
\end{pspicture}
}

\newcommand{\Gtwop}{
\begin{pspicture}[.4](-.5,0)(2,1)
\pscircle(0,.5){.35}
\pscircle(1.5,.5){.35}
\psline(.35,.5)(1.15,.5)
\put(-.7,.5){$i$}
\put(2,.5){$j$}
\put(.6,.7){$k$}
\end{pspicture}
}

\newcommand{\Ccurves}{
\psset{unit=1cm}
\begin{pspicture}[.4](-.6,-.8)(2.6,.8)
\pscircle(0,0){.35}\pscircle(2,0){.35}
\psline(-.05,-.05)(.05,.05)
\psline(-.05,.05)(.05,-.05)
\psline(1.95,-.05)(2.05,.05)
\psline(1.95,.05)(2.05,-.05)
\pscircle(0,0){.45}\pscircle(2,0){.45}
\pscircle(2,0){.55}
\psarc(0,0){.7}{45}{315}\psarc(2,0){.7}{225}{135}
\pscurve(.5,.49)(.6,.43)(.8,.4)(1.2,.4)(1.4,.43)(1.5,.49)
\pscurve(.5,-.49)(.6,-.43)(.8,-.4)(1.2,-.4)(1.4,-.43)(1.5,-.49)
\end{pspicture}
}
\newcommand{\Gtroisi}{
\psset{unit=.8cm}
\begin{pspicture}[.4](-.6,-2.8)(2.6,.8)
\psline(-.05,-.05)(.05,.05)
\psline(-.05,.05)(.05,-.05)
\psline(1.95,-.05)(2.05,.05)
\psline(1.95,.05)(2.05,-.05)
\pscircle(0,0){.6}
\pscircle(2,0){.6}
\pscircle(1,-1.73){.6}
\psline(.95,-1.78)(1.05,-1.68)
\psline(.95,-1.68)(1.05,-1.78)
\psline(1,-.577)(.51,-.29)
\psline(1,-.577)(1.49,-.29)
\psline(1,-.577)(1,-1.13)
\put(.1,-.1){${1}$}
\put(2.1,-.1){${2}$}
\put(1.1,-1.83){$3$}
\end{pspicture}
}
\newcommand{\Gtroisii}{
\psset{unit=.8cm}
\begin{pspicture}[.4](-.6,-2.8)(2.6,.8)
\put(.1,-.1){${1}$}
\put(2.1,-.1){${2}$}
\put(1.1,-1.83){$3$}
\put(1.6,-2.2){$A_{123}$}
\psline(-.05,-.05)(.05,.05)
\psline(-.05,.05)(.05,-.05)
\psline(1.95,-.05)(2.05,.05)
\psline(1.95,.05)(2.05,-.05)
\psarc(0,0){.6}{350}{310}
\psarc(2,0){.6}{230}{190}
\psarc(1,-1.73){.6}{110}{70}
\psline(.6,-.1)(1,-.33)
\psline(1.4,-.1)(1,-.33)
\psline(.385,-.475)(.8,-.7)
\psline(1.615,-.475)(1.2,-.7)
\psline(1.2,-.7)(1.2,-1.175)
\psline(.8,-.7)(.8,-1.175)
\psline(.95,-1.78)(1.05,-1.68)
\psline(.95,-1.68)(1.05,-1.78)
\end{pspicture}
}

\newcommand{\Gtrois}{
\psset{unit=.8cm}
\begin{pspicture}[.4](-.6,-2.8)(2.6,.8)
\put(.1,-.1){${1}$}
\put(2.1,-.1){${2}$}
\put(1.1,-1.83){$3$}
\put(1.6,-2.2){$A_{3}$}
\put(2.3,-.8){$A_{12}$}
\psline(-.05,-.05)(.05,.05)
\psline(-.05,.05)(.05,-.05)
\psline(1.95,-.05)(2.05,.05)
\psline(1.95,.05)(2.05,-.05)
\psline(.95,-1.78)(1.05,-1.68)
\psline(.95,-1.68)(1.05,-1.78)
\psarc(0,0){.6}{20}{340}
\psarc(2,0){.6}{200}{160}
\psline(.56,.2)(1.44,.2)
\psline(.56,-.2)(1.44,-.2)
\pscircle(1,-1.73){.6}
\end{pspicture}
}

\newcommand{\Ccalca}{
\begin{pspicture}[.4](-2,-1.5)(3.5,1.5)
\psline(-.1,-.1)(.1,.1)
\psline(-.1,.1)(.1,-.1)
\pscircle(0,0){1.5}
\psline(1.5,0)(3,0)
\put(-2,.1){$i$}
\put(2,.3){$k$}
\end{pspicture}
}
\newcommand{\Ccalcai}{
\begin{pspicture}[.4](-1.5,-1.5)(3.5,1.5)
\psline(-.1,-.1)(.1,.1)
\psline(-.1,.1)(.1,-.1)
\pscircle(0,0){1.5}
\psline(1.5,0)(3,0)
\put(-2,.1){$i'$}
\put(2,.3){$k$}
\end{pspicture}
}

\newcommand{\Ccalc}{
\begin{pspicture}[.4](-2,-1.5)(3.5,1.5)
\psline(-.1,-.1)(.1,.1)
\psline(-.1,.1)(.1,-.1)
\pscircle(0,0){1}
\put(-.6,.25){$\tilde p$}
\pscircle(0,0){1.5}
\psline(1.5,0)(3,0)
\put(-2,.1){$i$}
\put(2,.3){$k$}
\end{pspicture}
}

\newcommand{\Ccalci}{
\begin{pspicture}[.4](-1.5,-1.5)(3.5,1.5)
\psline(-.1,-.1)(.1,.1)
\psline(-.1,.1)(.1,-.1)
\psarc(0,0){1}{270}{90}
\put(.35,0){$\tilde p$}
\psarc(0,0){1.5}{270}{90}
\psarc(0,0){1.25}{118}{242}
\pscurve(0,1.5)(-.3,1.4)(-.6,1.1)
\pscurve(0,1)(-.3,1.02)(-.6,1.1)
\pscurve(0,-1.5)(-.3,-1.4)(-.6,-1.1)
\pscurve(0,-1)(-.3,-1.02)(-.6,-1.1)
\psline(1.5,0)(3,0)
\put(-1.8,.1){$i'$}
\put(2,.3){$k$}
\put(1.1,1.4){$i$}
\put(1.3,-1.3){$i$}
\end{pspicture}
}

\begin{document}
\title[Integral Bases and unimodular representations]{Integral bases for TQFT Modules and unimodular representations 
of mapping class groups}
 
\author{ Patrick M. Gilmer}
\address{Department of Mathematics\\
Louisiana State University\\
Baton Rouge, LA 70803\\
USA}
\email{gilmer@math.lsu.edu}

\author{ Gregor Masbaum}
\address{Institut de Math\'ematiques de Jussieu (UMR 7586 du CNRS)\\
Universit\'e Paris 7 (Denis Diderot) \\
Case 7012\\
2, place Jussieu\\
75251 Paris Cedex 05\\
FRANCE }
\email{masbaum@math.jussieu.fr}

\author{Paul van Wamelen}
\address{Department of Mathematics\\
Louisiana State University\\
Baton Rouge, LA 70803\\
USA}
\email{wamelen@math.lsu.edu}

\date{July 8, 2002}

\begin{abstract} We construct integral bases for the $SO(3)$-TQFT-modules 
of surfaces in genus one  and two at roots of unity of 
prime order and show that the corresponding mapping class group 
representations preserve a unimodular Hermitian form over a ring 
of algebraic integers. For higher genus surfaces the Hermitian form 
sometimes must be non-unimodular. In one such case, genus~$3$
and $p=5,$ we still give an explicit basis. 
\end{abstract}

\maketitle

\section{Introduction}

Integrality properties of Witten-Reshetikhin-Turaev quantum invariants 
of $3$-manifolds have been 
studied intensively in the last 
several
years. H. Murakami \cite{Mu1,Mu2} showed that 
the $SU(2)$- 
and $SO(3)$-invariants at a root of unity $q$ of prime order are algebraic 
integers. This was reproved in \cite{MR} and generalised to all classical Lie types in \cite{MW,TY} and then to all Lie types in \cite{Le}. These integrality properties are crucial for establishing the relationship of the invariants with the Casson invariant \cite{Mu1,Mu2} and with the perturbative invariants or Ohtsuki series 
\cite{OhCambridge,Oh1,Le}.

Quantum invariants fit into Topological Quantum Field Theories (TQFT). This means in particular that 
there are representations of mapping class groups associated
with them. 
(Actually the representations are usually only projective-linear; equivalently, one has to consider certain central extensions of mapping class groups here.) If a $3$-manifold $M$ is presented as a Heegaard splitting where two handlebodies are glued together by a diffeomorphism $\varphi$ along their boundary, the quantum invariant of $M$ can be recovered from the representation of $\varphi$ on the TQFT-vector space $V(\Si)$ associated to the boundary surface $\Si$. 

The TQFT-representations are finite-dimensional and can be defined over a finite extension of the cyclotomic number field $\mathbb{Q}(q)$, where the quantum parameter $q$ is a root of unity. They also preserve a non-degenerate Hermitian form $\langle \ , \ \rangle_\Si$ on $V(\Si)$ (which may or may not be unitary; this usually depends on the choice of the embedding of the cyclotomic field  into $\mathbb{C}$). 

A quite striking result was recently announced by Andersen  \cite{An} who proved that in the $SU(n)$ case the representations are asymptotically faithful (here asymptotically means letting the order of $q$ go to infinity). At a fixed root of unity, they are certainly not faithful, as  Dehn twists are always represented by matrices of finite order. Roberts \cite{R} showed that the representations are irreducible in the $SU(2)$-case if the order of $q$ is prime. 

An interesting question is to determine the image of the  mapping class group in the TQFT-representations.  For the $SU(2)$ and $SO(3)$-theories, the first author proved that in the genus one case the image is a finite group \cite{G1}. However in higher genus, the image is not finite \cite{Fu}; in fact, it contains elements of infinite order \cite{Madeira}.  One might hope that this image is equal to the linear transformations which are automorphisms of some (yet to be found) structure, just as a linear transformation of the homology group $H_1(\Si;\mathbb{Z})$ is represented by a mapping class if and only if it preserves the intersection form.

In this paper we are concerned with integrality properties of the TQFT representations. For simplicity, we restrict ourselves to the $SO(3)$ case; specifically, we use a variant of the $V_p$-theories of \cite{BHMV2} with $p$ an odd prime. Here, $p$ is the order of the root of unity $q$. Let $\BO$ denote the ring of algebraic integers in the cyclotomic ground field. The main idea to obtain an integral structure on the TQFT already appears in \cite{G}. Namely, we define an $\BO$-submodule $\BS(\Si)$ of the TQFT-vector space $\Vp(\Si)$ as the $\BO$-span of vectors represented by connected $3$-manifolds with boundary $\Si$. The point of this definition is that the submodule $\BS(\Si)$ is clearly preserved under the mapping class group. 

It was shown in \cite{G} that $\BS(\Si)$ is always a free finitely generated $\BO$-module. One can also rescale the Hermitian form $\langle \ , \ \rangle_\Si$ on $\Vp(\Si)$ to obtain a non-degenerate $\BO$-valued  form $(\ ,\ )_\Si$ on $\BS(\Si)$.  This relies on the integrality results for the \hbox{$3$-manifold} invariants mentioned above. 

The form $(\ ,\ )_\Si$ is again preserved by the mapping class group. In particular, the image of the mapping class group in the TQFT-representation $V_p(\Si)$ lies in the subgroup preserving a lattice defined over $\BO$ and a non-degenerate Hermitian form on it. 

In what sense is $\BS$ a TQFT defined over $\BO$? For instance one might hope that the form $(\ ,\ )_\Si$ was unimodular. Here, we show that this is indeed the case in genus one and two. This is a consequence of our main result which is to describe explicit bases of 
$\BS(\Si)$
in genus one and two. 

In fact, we will describe two quite different bases in genus one. The first basis is  given in Theorem~\ref {firstbasis}. It is $$\{ \omega, t(\omega) , t^2(\omega) , \ldots , t^{d-1}(\omega)\}$$  where $\omega$ is the element appearing in the surgery axiom of the $\Vp$-theory and $t$ is the twist map. It is easy to see that these elements lie in $\BS(S^1\times S^1)$ and we use a Vandermonde matrix argument to show that they form a basis. A crucial step is to show that the Hermitian form $(\ ,\ )_{S^1\times S^1}$ is unimodular with respect to this basis. 

The second basis given in Theorem~\ref{secondbasis} is of a quite different 
nature.  
It is $$\{ 1, v ,v^2 , \ldots , v^{d-1}\}$$ where $v=(z+2)/(1+A)$ (here $z$ is represented by the core of the solid torus). We call it the $v$-basis. This time,  it is not even obvious {\em a priori} that its elements lie in $\BS(S^1\times S^1)$. The proof involves two different arguments. One is to show that the $\BO$-span of the $v$-basis is stable under the twist map. This is shown in Section~\ref{Kauf}. In fact, we prove it in the more general context where the skein variable $A$ is an indeterminate rather than a root of unity.
The second ingredient is to express $\omega$ in the $v$-basis and thereby relate the $v$-basis to the first basis. This is done in Section~\ref{2nd}.

The $v$-basis lends itself nicely to finding bases in higher genus. In Section~\ref{g2}, we describe a basis of $\BS(\Sigma)$ in genus two consisting of $v$-colored links in a genus two handlebody. These links are described by arrangements of curves in a twice punctured disk. Again, the unimodularity of the Hermitian form with respect to this basis is a crucial step in the argument.

In principle this method can be used to study $\BS(\Sigma)$ in higher 
genus as well. It turns out, however, that  the Hermitian 
 form $(\ ,\ )_{\Si}$ is 
not always unimodular.
For example, a simple argument given in Section \ref{non-uni} shows that it cannot be unimodular for surfaces of genus 3  and 5,  assuming $p \equiv 5 \pmod{8}.$ 

In this paper we 
will not
attempt to deal with the higher genus case in general. We only give in Section~\ref{g3}   a basis of $\BS(\Sigma)$ for a surface of genus three when $p=5.$ Although in this case the Hermitian form is not 
unimodular, it is nearly so. This allows us to find a basis  easily
in this one case.

Note that our definition of $\BS(\Si)$ is analogous to the construction of integral modular categories in \cite{MW}; in both cases one constructs integral structures by considering the span, over the subring of algebraic integers of the coefficient field, of the morphisms of the geometrically defined  category (tangles in the case of \cite{MW}, $3$-dimensional cobordisms  in the case at hand). It might be that this is not always enough: It is conceivable that one might be able to enlarge $\BS(\Si)$ in some 
way
to make the form always unimodular; however this enlargement would not be generated by $3$-cobordisms anymore.

We conclude the paper by showing how the $\BS$-theory defined over $\BO$ can be used to prove a divisibility result for the Kauffman bracket of
links in $S^3$. This generalizes a result of Cochran and Melvin \cite{CM} for zero framed links 
(see also \cite{OhCambridge, KS}).

\vskip 8pt
\noindent {\em Notational conventions.} Throughout the paper, $p\geq 3$ will be an odd integer, and we put $d=(p-1)/2$. From Section~\ref{elem} onwards, $p$ is supposed to be prime.

\section{The twist map on the Kauffman Bracket module of a solid torus}\label{Kauf}
In this section we define a sequence of submodules $K(n)$ of the Kauffman Bracket skein module of the solid torus $S^1 \times D^2$ and show that they are preserved under the twist map. We use the notations of \cite{BHMV1}. 

Suppose $R$ is a commutative ring with identity and an invertible
element $A.$ The universal example is $R=\BZ[A,A^{-1}]$ which we also
denote by $\BZ[A^{\pm}]$.  Recall that the Kauffman bracket skein
module $K(M,R)$ of a $3$-manifold $M$ is the free $R$ module generated
by isotopy classes of banded links in $M$ modulo the submodule
generated by the Kauffman relations.
  
 We let $z$ denote the skein element of $K(S^1 \times D^2,R)$ given by
the banded link $S^1\times J,$ where $J$ is
a small
arc in
 the interior of  
$D^2.$

 As is well known, $K(S^1 \times D^2,R)$ is  a free $R$-module on the 
nonnegative powers of $z,$ where $z^n$ means $n$ parallel copies of $z$. 
This also makes 
$K(S^1 \times D^2,R)$ into an $R$-algebra isomorphic to the polynomial ring
$ R [z].$ 
 
Let $t:K (S^1 \times D^2,R) \rightarrow  K (S^1 \times D^2,R)$  denote the twist map induced by a full right handed twist on the solid torus. It is well known (see {\em e.g.} \cite{BHMV1}) that there is a basis $ \{e_i\}_{i\ge0}$ of 
eigenvectors for the twist map. It is defined recursively by
 \begin{equation}\label{ei}
e_0=1, \quad e_1=z,  \quad e_i= z e_{i-1}- e_{i-2} ~.
\end{equation}
 The eigenvalues are given by  
\begin{equation}\label{mui}t(e_i)=\mu_i e_i, \text{ where } \mu_i=(-1)^i A^{i^2+2i} ~.\end{equation} 

\begin{de} {\em Let ${K}(n)$ denote the $\BZ[A^{\pm}]$-submodule of $K  (S^1 \times D^2, \BZ[A^{\pm},\frac 1 {1+A}])$ generated by $\{1,v,v^2,\ldots,v^n\}$, where $$v= \frac {z+2}{1+A}~.$$
}  
\end{de}

\begin{thm} \label{tK} The twist map $t$  sends 
$K(n)$ to itself.
\end{thm}
\begin{proof} Consider the basis $ \{ (z+2)^i\}_{i\ge0}$ of 
$K (S^1 \times D^2,\BZ[A^{\pm}]).$
The following Lemma gives the change of basis formulas.

\begin{lem}
\label{changeofbasis} 
For each $n \ge 1,$
\[(z+2)^{n-1}= \sum_{k=1}^n \binom{2n}{n-k} \frac k n e_{k-1}\]
\[ e_{n-1}= \sum_{i=1}^n (-1)^{n-i}  \binom{n+i-1}{n-i} (z+2)^{i-1}\]
\end{lem}

\begin{proof} Prove each separately by
induction on $n$ using the recursion formula (\ref{ei}).

\end{proof}

\begin{rem} {\em It follows that $\binom{2n}{n-k} \frac k n \in \BZ,$ which can also be seen directly:}
\[ \binom{2n}{n-k} \frac k n=
\binom{2n}{n-k} (1-\frac {n-k} n) =
\binom{2n}{n-k}- 2\binom{2n-1}{n-k-1}~.\]
\end{rem}

It is enough to show Theorem~\ref{tK} for the endomorphism $-At$ in place of $t$. Let us compute $-At$ in the basis $(z+2)^n.$ Note that $-At(e_{i-1})=(-A)^{i^2} e_{i-1}$.
\begin{align*}
-At\left((z+2)^{n-1}\right) & = \sum_{k=1}^n \binom{2n}{n-k} \frac k n  (-A) t(e_{k-1}), \\
 &= \sum_{k=1}^n \binom{2n}{n-k} \frac k n (-A)^{k^2} \sum_{i=1}^k (-1)^{k-i} \binom {k+i-1}{k-i} (z+2)^{i-1},\\
&= \sum_{i=1}^n  (-1)^i 
\left( \frac 1 n
\sum_{k=i}^n k \binom {2n}{n-k} \binom {k+i-1}{k-i} A^{k^2} \right)(z+2)^{i-1},\\
&=  \sum_{i=1}^n (-1)^{i} S_{1,i,n}(A) (z+2)^{i-1}, 
\end{align*}

Here, for $ m \ge 1,$ we define  
\[S_{m,i,n}(A)= \frac 1 n  \sum_{k=i}^n  k^m \binom {2n} {n-k}   \binom {k+i-1}{k-i} A^{k^2}\in \BZ[A].\] 

\begin{lem} \label{at  -1}
\begin{equation*} 
S_{1,i,n}(-1) = 
\begin{cases}
(-1)^n, &i=n\\
0  &i \ne n
\end{cases}
\end{equation*}
\end{lem}

\begin{proof} 
If we put $A=-1,$ then all $\mu_i=1$ and hence $-At$ is the identity.
\end{proof}

The following  formula is a very special
case of a transformation formula for terminating hypergeometric series
due to Bailey \cite[Formula 4.3.1]{B}. This was pointed out to us
by Krattenthaler's {\tt HYP} package, \cite{hyp}.

\begin{lem} \label{recurs}
\begin{equation*} 
S_{m,i,n}= i^2 S_{m-2,i,n} + 2i (2i+1) S_{m-2,i+1,n}
\end{equation*}
\end{lem}

\begin{proof} 
\begin{align*}
 S_{m,i,n} 
&=\frac 1 n  \sum_{k=i}^n  k^{m-2} \left( (k+i) (k-i) +i^2 \right) \binom {2n} {n-k}   \binom {k+i-1}{k-i} A^{k^2} \\
&= i^2 S_{m-2,i,n}  + \frac 1 n \sum_{k=i+1}^n  k^{m-2}  (k+i) (k-i)  \binom {2n} {n-k}   \binom {k+i-1}{k-i} A^{k^2} \\
& \text{ (the term with $k=i$ is zero) }\\
&= i^2 S_{m-2,i,n}  +  \frac 1 n 2 i (2i+1)\sum_{k=i+1}^n  k^{m-2} \binom {2n} {n-k}   \binom {k+i}{k-i-1} A^{k^2} \\
&= i^2 S_{m-2,i,n} + 2i (2i+1) S_{m-2,i+1,n}
\end{align*}

Here we use the simple identity:
\[ (k+i)(k-i) \binom {k+i-1}{k-i} = 2i (2i+1) \binom {k+i}{k-i-1} \]
\end{proof}

\begin{prop}\label{eq2.7} \label{divis} $S_{1,i,n}(A)$ is divisible by $(1+A)^{n-i}$ in $\BZ[A]$ for $i \le n.$
\end{prop}

\begin{proof} It suffices to show:
\begin{equation}\label{eq2}  
\left[ \left( \frac {d} {dA} \right) 
^{k}
S_{1,i,n}(A) \right]_{A=-1} =0
\end{equation}
for all $k=0,1,\ldots , n-i-1$.
Note that  
$\frac {d} {dA}  S_{m,i,n}= A^{-1} S_{m+2,i,n} \in \BZ[A].$ Thus

\[
 \left( \frac {d} {dA} \right) 
^k
S_{1,i,n} \in \text{ Span}_{\BZ[A^{\pm}]}\{ S_{m,i,n} \,|\, m \text{  odd}, \, 1 \le m \le 2n-2i-1\} 
\]
for all $k$ in the required range. 
\noindent But using Lemma~\ref{recurs} one may decrease $m$ at the cost of increasing $i,$ and see that 
\[
 \text{ Span}_{\BZ[A^{\pm}]}\{ S_{m,i,n} \,| \, m \text{  odd},\, 1 \le m \le 2n-2i-1\} \subseteq  \text{ Span}_{\BZ[A^{\pm}]}\{ S_{1,j,n} \,|\, j < n\}.
\]

\noindent  By Lemma~\ref{at -1},  $S_{1,j,n}(-1) =0,$ for  $ j < n,$ and (\ref{eq2})  follows.
\end{proof}

\noindent{\em Proof of Theorem~\ref{tK}.} We have
$$t(v^{n-1}) = - A^{-1}  \sum_{i=1}^n s_{i,n}(A) v^{i-1},$$ where 
$s_{i,n}(A)= (-1)^i (1+A)^{i-n}S_{1,i,n}(A)$  lies in $\BZ[A]$ by Proposition~\ref{eq2.7}. 
\end{proof}

\begin{rem}{\em Theorem~\ref{tK}  remains valid if we take $A$ to be a root of unity, other than $-1,$ rather than an indeterminant. }
\end{rem}

\begin{rem}{\em Let $\tilde {K}(n)$ be defined as $K(n)$ but with 
$v=({z+2})/({1+A})$ replaced with $\tilde v= ({z+2})/({1-A^2})$. 
Then a similar argument shows that $\tilde {K}(n)$ is stable under 
$t^2$, the square of the twist map. (To see this, one should replace
 $-At$ with $A^2 t^2$ in the above
and express everything in terms of $q=A^2$.
This  leads to polynomials 
       $\tilde S_{m,i,n}(q)$
defined similarly as the $ S_{m,i,n}(A)$ except
 that $A$  is replaced with 
$q$
and an extra factor of $(-1)^k$ is 
inserted in the sum. The remainder of the argument is the same.)
}
\end{rem}

\section{The $SO(3)$-TQFTs}

Let $p\geq 3$ be an odd integer. (In this section, $p$ need not be prime.) We consider a variation of the $2+1$ dimensional cobordism category considered in \cite{BHMV2} whose objects are closed oriented surfaces (with extra structure) with a (possibly empty) collection of banded points ($=$ small oriented arcs) colored by integers in the range $[0,p-2].$ The morphisms are (equivalence classes of) oriented 3-dimensional manifolds (with extra structure)  with $p$-admissibly colored banded trivalent graphs. (Two morphisms are considered equivalent if they are related by a homeomorphism respecting the boundary identifications.) For the definition of $p$-admissibility in the $p$-odd case see \cite[Theorem 1.15]{BHMV2}; see also Section~\ref{g2}.

The variation consists of replacing the $p_1$-structures of \cite{BHMV2} with structures put forward by Walker \cite{W} and Turaev \cite{Tu}.
Surfaces are equipped with a Lagrangian subspace of their first homology. We use homology with  rational coefficients when considering 
Lagrangian subspaces. 
Cobordisms
are equipped with integer weights, as well as Lagrangian subspaces for the target and source. This is also described in \cite{G}. We will denote this category by 
$\mathcal{C}.$ We call the objects of this category e-surfaces, and call  the morphisms 3-e-manifolds.

The procedure of   \cite{BHMV2} defines a TQFT-functor  $\Vp$ 
on $\mathcal{C}$ over
a commutative ring $R$ containing 
$p^{-1}$, 
a primitive 2pth root of unity $A$ and a 
solution of 
        $\kappa^2= A^{ -6- p(p+1)/2 }.$
The number $\kappa$ here plays the role of $\kappa^3$ in \cite{BHMV2}. Here we use the term TQFT slightly loosely as the tensor product axiom does not hold
unless only even colors are used in the cobordism category. The even colors correspond to 
irreducible
representations of $SU(2)$ which lift to $SO(3)$. Therefore the $\Vp$-theory for odd $p$ is considered a $SO(3)$ variant of the Witten-Reshetikhin-Turaev $SU(2)$-TQFT.  

For us it is convenient to use odd colors as well as  even colors. However, if we insist that only even colors  be used in coloring the banded points on the surfaces,  
then we do obtain  
an 
honest TQFT with the tensor product axiom, but we are still
allowed us to use the language of odd colors to describe states. 
This will be 
useful in 
Sections~\ref{g2} and~\ref{g3}.

If $M$ is a 3-e-manifold viewed as morphism from $\Si$ to $\Si'$ in $\mathcal{C}$, we denote the associated endomorphism from $\Vp(\Si)$ to $\Vp(\Si')$  by 
$\Zp(M)$. (It is denoted by $(Z_p)_M$ in \cite{BHMV2}).
If $M$ is a closed 3-e-manifold viewed as morphism from $\emptyset$ to $\emptyset,$ $\Zp(M)$ induces multiplication by a scalar from $R=  V(\emptyset).$ 
This scalar is denoted by $\langle M\rangle.$  If $M$ is a  3-e-manifold viewed as morphism from $\emptyset$ to $\Sigma,$  let $[M]$ denote  $\Zp(M)(1) \in V(\Sigma).$  ($[M]$ is denoted by $Z_p(M)$ in \cite{BHMV2}).  We call such an element $[M]$ a {\em  vacuum state.} If $M$ is connected, $[M]$ is called a  connected vacuum state.

The modules $\Vp(\Si)$ are always free over $R$. They also carry a nonsingular Hermitian form \cite{BHMV2} :
\[\langle \ ,\ \rangle _{\Sigma}: \Vp(\Sigma) \times \Vp(\Sigma) \rightarrow R\]
given by 
\[\langle [N_1],[N_2]\rangle _{\Sigma}= \langle N_1 \cup_{\Sigma} -N_2\rangle~.\]
Here $-N_2$ is the 3-e-manifold obtained by reversing the orientation, multiplying the weight by $-1, $ and leaving the Lagrangian on the boundary alone.

If $\Sigma$ is an e-surface with no colored points, and $H$ is a handlebody
(weighted zero)
with boundary $\Sigma,$ then $\Vp(\Sigma)$ has a specified isomorphism to a quotient of the skein module $K(H,R)$
\cite[p. 891]{BHMV2}. In fact if $H$ is a subset of $S^3$ then two skeins represent  the same element
if 
and only if 
they are equal as ``maps of outsides'' in Lickorish's phrase
\cite{L}. 

Let $S^1 \times S^1$ denote an e-surface of genus one with no colored points.
Let $d$ denote $(p-1)/ 2.$ It turns out that $d$ is the dimension or rank of 
$\Vp(S^1 \times S^1).$ 
In fact, the module  $\Vp(S^1 \times S^1) $ is isomorphic as an $R$-module to the quotient of $K(S^1 \times D^2,R)=R[z]$ by the ideal generated by $e_d- e_{d-1}\in R[z].$  It follows \cite[p.696]{BHMV1} that $e_{p-1}=0$ in   $\Vp(S^1 \times S^1) $ and $e_{d+i}= e_{d-1-i}.$ Thus the module $\Vp(S^1 \times S^1) $ has indeed rank $d$ with the  basis $\{e_0,e_1,\ldots ,e_{d-1}\}$. Note that this basis is the same, up to reordering, as the even basis $\{e_0,e_2,\ldots, e_{p-3}\}.$

Let $\Sigma_g$ denote an e-surface of genus $g$ with no colored points on the boundary. 
The rank of the free module $\Vp(\Sigma_g)$ is given by the formula
\cite[Cor. 1.16]{BHMV2}
\[ \text{rank}\left(\Vp(\Sigma_g)\right) = \left( \frac p 4\right)^{g-1} \sum_{j=1}^{d}
 \left( \sin \frac {2 \pi j}{p} \right)^{2-2g} .\] This is the same as $2^{-g}$ times the dimension of $V_{2p}(\Si_g)$ (this fact comes from a tensor product formula, see 
         \cite[Thm. 1.5]{BHMV2}).
 Note that $V_{2p}(\Si_g)$ is an $SU(2)$-TQFT module, with dimension given by the $SU(2)$ Verlinde formula at level $p-2$ (where the colors are again the set of integers in the range $[0,p-2]$).

In genus 2, we have
\[ \text{rank}\left(\Vp(\Sigma_2)\right)= \frac {d (d+1)(2d+1)}{6} \]
as will be seen by an explicit counting argument in Section~\ref{g2}.

\section{Some facts from elementary number theory}\label{elem}

In the remainder of this paper,  we assume $p$ is an odd prime. We continue to use the notation $d=(p-1)/2$.

In this section, we collect some notation and a few elementary number-theoretical facts. To be specific we pick particular values for $A$ and $\kappa$. 
We put 
$A= \zeta_{2p}$ where $\zeta_n = e^{2 \pi i/{n}},$ and also use 
the 
notation\footnote{Warning: In many places ({\em e.g.} in 
\cite{MR}), $q$ denotes $A^4$ rather than $A^2.$}
 $$q=A^2.$$ We may then take $\kappa = A^{-3}(-i)^{({p+1})/{2}}.$ 
Note that $\zeta_{2p}\in \BZ[\zeta_p]$.
Thus the coefficient ring is  
$R =\BZ[\zeta_p, \frac 1 p]$ if $p \equiv -1 \pmod{4},$ and 
$R =\BZ[\zeta_{p},i, \frac 1 p]=\BZ[\zeta_{4p}, \frac 1 p]$ if $p \equiv 1 \pmod{4}.$ Of course, the coefficient ring remains  unchanged if $A$ is replaced by another primitive $2p$-th root of unity, and $\kappa$ is changed accordingly.

We let $\eta$ denote $\langle S^3\rangle$, the invariant of $S^3$ with weight zero, and put 
$\mathcal{D}= 
\eta^{-1}$.
Then 
using equations on \cite[p.897]{BHMV2}
\begin{equation}\label{D}
  \mathcal{D} = \frac{i\ \sqrt{p} }{q-q^{-1}} = \frac{i^{\frac {p+1} 2} }{q-q^{-1}} 
\left(\frac 1 2 \sum _{m=1}^{2p}(-1)^m A^{m^2}\right)
\end{equation}
In particular
\begin{equation}\label{DD} \mathcal{D}^2= \frac {-p}{ ({q-q^{-1}})^2}~.
\end{equation}

 We denote by $\mathcal{O} $ the  ring of integers in $R.$ Note that $\mathcal{O}=\BZ[\zeta_p]$ if $p \equiv -1 \pmod{4},$ and 
$\mathcal{O}  =\BZ[\zeta_{p},i]=\BZ[\zeta_{4p}]$ if $p \equiv 1 \pmod{4}.$
 
The following notation will be useful.
If $x,y$ are elements of $\BO$ (or, more generally, of its quotient field), 
we write $x\sim y$ if there exists 
a unit 
$u\in \BO$ 
such that $x=uy$. 

\begin{lem}\label{num}
  \begin{itemize}
\item[(i)] $1-A$ is a unit in $\BO$, and $1-q\sim 1+A$. 
\item[(ii)] 
One has $\mathcal{D}\in \BO$. Moreover, 
$\mathcal{D}\sim (1-q)^{({p-3})/2}=(1-q)^{d-1}$.
\item[(iii)] The quantum integers $[n]=(q^n-q^{-n})/(q-q^{-1})$ are units for $1\leq n\leq p-1$.
\item[(iv)] If $0\leq i,j \leq d-1$ and $i\neq j$, then the twist coefficients $\mu_i$ (see (\ref{mui})) satisfy $\mu_i-\mu_j\sim 1-q$.
 \item[(v)] Put $\lambda_i=-q^{i+1}-q^{-i-1}$. If $0\leq i \leq d-1$, then $\lambda_0-\lambda_i\sim (1-q)^2$.
  \end{itemize}
\end{lem}
\begin{proof}  The fact that $1-A$ is a unit follows easily from the fact that $A$ is a zero of the $2p$-th cyclotomic polynomial $1-X+X^2-\ldots +X^{p-1}$. This proves (i).  It is well-known that  $p\sim (1-q)^{p-1}$ (see {\em e.g.} \cite[Lemma 3.1]{MR}). Together with 
Formulas (\ref{D})
and
 (\ref{DD}), this shows (ii). Observing that $[n]\sim
1+q^2+\ldots +q^{2n-2}$,
(iii) is also shown in \cite[Lemma 3.1]{MR}.  For (iv), observe that $\mu_i=\mu_{p-2-i}$ so that the set of $\mu_i$ in question is equal to the set of $\mu_{2i}=q^{
2 i^2+2i}$ for $i=0,1,\ldots d-1$. These powers of $q$ are all distinct, which implies $$\mu_i-\mu_j \sim 1-q^n \sim 1-q$$ for some $0<n<p$. This proves (iv). Similarly, (v) follows from $\lambda_0-\lambda_i\sim(1-q^{i+2})(1-q^i)$. 
\end{proof}
\begin{rem}{\em It is well-known that $1-q$ is a prime in $\BZ[q]=\BZ[\zeta_{p}]$. But if $p \equiv 1 \pmod{4},$ then $1-q$ is not a prime in $\BO=\BZ[\zeta_{4p}]$ (it splits as a product of two conjugate prime ideals). 
}\end{rem}

\section{Associated Integral Cobordism Functors}

In \cite{G}, a cobordism functor from a restricted cobordism category to the category of free finitely generated $\mathcal{O}$-modules is described. Let $\mathcal{C''}$ denote the subcategory of  $\mathcal{C}$
defined by considering only nonempty connected surfaces and connected morphisms between such surfaces.
This represents a further restriction of $\mathcal{C}$ than that considered in \cite{G}, but it suffices for our purposes.

\begin{de}[\cite{G}] {\em If $\Sigma$ is a connected  e-surface, define $\BS(\Sigma)$ to 
be the $\mathcal{O}$-submodule of $\Vp(\Sigma)$ generated by connected vacuum states. If $N: \Sigma \rightarrow \Sigma'$ is a morphism
of $\mathcal{C''}$ then $\Zp(N)$ sends $ [M] \in \BS(\Sigma)$ to $[M \cup_{\Sigma}N] \in \BS(\Sigma').$ In this way we get a functor   
from $\mathcal{C''}$ to the category of
finitely generated  $\mathcal{O}$-modules.
We   also rescale the Hermitian form on 
$\Vp(\Sigma)$
to obtain an $\mathcal{O}$-valued  Hermitian form 
\[ (\ ,\ )_{\Sigma}: \BS(\Sigma) 
\otimes_\BO
\BS(\Sigma) \rightarrow \mathcal{O},\]
 defined by
\[ ([N_1],[N_2])_{\Sigma}= 
\mathcal{D} \langle [N_1],[N_2]\rangle_{\Sigma} =
\mathcal{D} \langle N_1 \cup_{\Sigma} -N_2\rangle
.\]
}\em
\end{de}
This form takes values in $\mathcal{O}$ by the integrality result for closed 3-e-manifolds \cite{Mu2,MR}. These theorems are also used in proving that
$\BS(\Sigma)$ is finitely generated \cite{G}. 

\begin{rem} {\em Over a Dedekind domain such as $\BO$, a finitely generated
torsion-free
 module is always projective, but it need not be free. (The typical examples are non-principal ideals in $\BO$.) Somewhat surprisingly, 
however, it turns out 
that
the modules $\BS(\Sigma)$ are always free. This is proved in \cite{G}. We will not actually make use of this fact in genus $1$ and $2$: freeness will follow from the construction of explicit bases.
} 
\end{rem} 
 \begin{de} {\em A Hermitian form on a projective $\BO$-module $S$ is called {\em non-degenerate} (or {\em non-singular}) if its  adjoint map $S\rightarrow S^*$ is injective. It is called {\em unimodular} if the adjoint map is an isomorphism.
}\end{de}

Note that if $S$ is free and $M$ is the matrix of the Hermitian form in some basis, then the form is non-degenerate (resp. unimodular) if $\det M$ is non-zero (resp. a unit in 
$\BO$).

In our situation, the form $ (\ ,\ )_{\Sigma}$ is always non-degenerate (since the original form $\langle \ ,\ \rangle_\Si$ on 
$\Vp(\Si)$
is). We will show that $ (\ ,\ )_{\Sigma}$ is unimodular in  genus $1$ and 
$2$.
\vskip 8pt
There is a standard basis $\{u_\sigma\}$  of $\Vp(\Sigma_g)$ given by p-admissible { \it even} colorings $\sigma$ of the graph 

\begin{center}
\Gn
\end{center}

\noindent 
(where there are $g$ loops)
embedded in a 3-e-handlebody $H_g$ of genus $g$
with boundary the e-surface $\Si_g$
 (see \cite[4.11]{BHMV2}).  
One may actually use any trivalent graph in $H_g$ to which $H_g$ deformation retracts. 
 (In the case $g=1$, 
this is
the same as the basis given by the elements $e_i$.) These basis elements lie 
in $\BS(\Sigma_g)$ because the denominators appearing in the Jones-Wenzl idempotents needed to expand colored graphs into skein elements are invertible in $\BO$ (see \cite{MR}). Warning: the $u_\sigma$   do {\em not}  generate $\BS(\Sigma_g)$ over $\BO$.  

\begin{prop}\label{detb} The elements $u_\sigma$ are orthogonal for the form $(\ ,\ )_{\Sigma_g}.$ Moreover, one has
$$
(u_\sigma,u_\sigma)_{\Sigma_g}\sim \BD^g\sim (1-q)^{(d-1)g}.
$$
\end{prop}
\begin{proof} By \cite[Theorem 4.11]{BHMV2}
one has that $\langle u_\sigma,u_\sigma\rangle_{\Sigma_g}$ is equal to 
$\eta^{1-g}=\BD^{g-1}$ times a product of non-zero quantum integers or 
their inverses,
which are units in $\BZ[q]$ by  Lemma~\ref{num}.
Since  the form $(\ ,\ )_{\Sigma_g} $ is just a rescaling of the form 
$\langle \ ,\ \rangle_{\Sigma_g} $, the result follows.
\end{proof}

One of the reasons to study the  form $(\ ,\ )_{\Sigma_g} $ is that it is preserved by the TQFT-action of the mapping class group. More precisely, let  $\tilde{\Gamma}(\Sigma)$ denote the central extension of the mapping class group $\Gamma(\Sigma)$  of $\Sigma$  realized by the subcategory of  
$C''$
consisting 
of
e-manifolds homeomorphic to  $\Sigma \times I$ such that the colored graph is given by $I$ times the
colored banded points of $\Sigma.$ This homeomorphism need not respect the boundary identification at 
$\Sigma \times \{1\},$ but should respect the boundary identification at $\Sigma \times \{0\}.$
In fact considering this boundary identification  at $\Sigma \times \{1\},$  defines  the quotient homomorphism from $\tilde{\Gamma}(\Sigma)$ to $\Gamma(\Sigma),$ which has kernel $\BZ$ given by the integral weights on $\Sigma \times I$ with standard boundary identifications. The group $\tilde \Gamma(\Sigma)$  is isomorphic to the signature extension (see {\em e.g.} Atiyah \cite{At}, Turaev \cite{Tu}.) This extension  can be described nicely using skein theory \cite{MR1}. 

\begin{prop}  The group $\tilde \Gamma(\Sigma)$ acts on $\Vp(\Sigma)$ preserving the $\mathcal{O}$-lattice $\BS(\Sigma)$ and the $\BO$-valued Hermitian form 
$ (\ ,\ )_{\Sigma}$.
\end{prop}
\begin{proof} This follows from the definition of $\BS(\Sigma)$ and the fact that the group $\tilde \Gamma(\Sigma)$ preserves the original Hermitian form $ \langle\ ,\ \rangle_{\Sigma}$.
\end{proof}

The module $\BS(\Sigma)$ can be described using the notion of `mixed graph'. Recall the element $$\omega= \BD^{-1} \sum_{i=0}^{d-1}  \langle e_i\rangle e_i \in K(S^1 \times D^2,R).$$  Here $\langle e_i\rangle=(-1)^i [i+1]$. It plays an important role in the surgery axioms of the $\Vp$-theory. 

By a {\em mixed graph} in a weighted 3-manifold $M, $ we mean a trivalent banded graph in $M$ whose simple closed curve components may possibly be colored $\omega$ or by integer colors in the range $[ 0,p-2]$ and whose other edges are colored p-admissibly by integers in the range $[ 0,p-2]$. 
A mixed graph can be expanded multilinearly into a $R$-linear combination of colored graphs. The result should be thought of as a superposition of e-morphisms.   If the graph is a link and every component is colored $\omega,$ we say the link is $\omega$-colored. A mixed graph in a handlebody $H$ specifies an element in $\Vp(\partial H).$

\begin{thm}\label{ml} A mixed graph in a 
connected
 3-e-manifold $M$ with boundary $\Sigma$ represents an element of $\BS(\Sigma).$ If $H$ is a 
3-e-handlebody 
with boundary
e-surface 
$\Sigma$ then $\BS(\Sigma)$  is generated over $\mathcal{O}$ by elements specified by mixed graphs in $H.$ 
\end{thm}

\begin{proof} The first statement follows from  the fact that  $\Vp$ satisfies the surgery axiom (S2) \cite[p 889]{BHMV2}. The second statement follows from the fact that any connected 3-manifold with boundary $\Sigma$ can be obtained by a sequence of 2-surgeries to  $H$  \cite[Proof of Lemma p. 891]{BHMV2}.
\end{proof}

\begin{rem} \label{5.7} {\em Suppose that we know that some collection 
$T$ of elements of $\Vp(\Sigma)$ lie in the $\mathcal{O}$-lattice  
$\BS(\Sigma).$ Then 
$\text{Span}_\mathcal{O}(T)$ is a $\mathcal{O}$-sublattice of
$\BS(\Sigma).$ 
This sublattice
might not be invariant under
$\tilde{\Gamma}(\Sigma)$.  Let $G=\{g_i\} \in \tilde \Gamma(\Sigma)$
be a finite set of elements whose image in $\Gamma(\Sigma)$ generate.
The sequence of submodules of $\BS(\Sigma)$:
$\text{Span}_\mathcal{O}(T),$ $\text{Span}_\mathcal{O}(T\cup G(T)),$
$\text{Span}_\mathcal{O}(T\cup G(T) \cup G(G(T))),$ , $ \ldots $ etc.
must stabilize in an $\mathcal{O}$-sublattice of $\BS(\Sigma)$ which
is invariant under the mapping class group. This procedure is well
suited to computer investigation.  The basis given in
Section~\ref{2nd} was originally found by this procedure. We used the
computer program Kant \cite{D} starting with $T= \{e_0,e_1,\ldots,
e_{d-1}\}$ in 
$\BS(S^1\times S^1).$ 
}\end{rem}

\section{First integral basis in genus 1} \label{1st}

By a slight abuse of notation, we let $\omega$ denote the element in $\BS(S^1 \times S^1)$ given by coloring the core of $S^1 \times D^2$ with $\omega.$
Let $t$ also denote the induced map on $\Vp(S^1 \times S^1)$ given by giving $S^1 \times D^2$ a full right 
handed
twist. Note that $t^n(\omega)\in \BS(S^1 \times S^1)$ for all $n$.

\begin{thm} \label{firstbasis} $\{ \omega, t(\omega) , t^2(\omega) , \ldots , t^{d-1}(\omega)\}$ is a basis for the  module $\BS(S^1 \times S^1).$ The form $( \ ,\ )_{S^1 \times S^1}$ is unimodular.
\end{thm}

Note that it follows in particular that the $\BO$-span of $\{ \omega, t(\omega) , t^2(\omega) , \ldots , t^{d-1}(\omega)\}$ is stable under the action of the mapping class group $\tilde{\Gamma}(
S^1 \times S^1)
$. 

\begin{proof}
Recall that $\mu_i=(-1)^{i} A^{i^2+2i}$ denotes the eigenvalue of $e_i$ under the twist map  $t$.
We have that
\[ t^j (\omega)= \BD^{-1} \sum_{i=0}^{d-1} \langle e_i\rangle \mu_i^j e_i~.\] 
Note that $\langle e_i\rangle=(-1)^i [i+1]$ is a unit by Lemma~\ref{num}(iii). The matrix  $W$ which expresses $\{ \omega, t(\omega) , t^2(\omega) , \ldots , t^{d-1}(\omega)\}$ in terms of   $\{e_0,e_1,\ldots e_{d-1}\}$  has as determinant a unit (the product of the $\langle e_i\rangle$) times $\BD^{-d}$  times the determinant of the 
Vandermonde matrix $[  \mu_i^j ]$ where $0 \le i,  j\le d-1.$ 
Moreover by Lemma~\ref{num}(iv)
\[\det [  \mu_i^j ] = \pm \prod_{i<j} (\mu_i-\mu_j)\sim  (1-q)^{{d(d-1)}/2}~.\] 

As $\BD \sim   (1-q)^{d-1},$ we conclude that 
\begin{equation} \label{detW} \det W\sim (1-q)^{- {d(d-1)}/2}~.
\end{equation} 

In particular, this determinant is non-zero, hence the $t^j(\omega)$ are linearly independent. Let $\mathcal {W}$ denote the $\mathcal{O}$-module spanned by the $t^j(\omega)$. Clearly $\mathcal W \subset \BS(S^1 \times S^1)$.  
Now we know by Proposition~\ref{detb} that $(e_i,e_i)\sim (1-q)^{d-1}$ (here we simply write $(\ ,\ )$ for the 
Hermitian
form $(\ ,\ )_{S^1 \times S^1} $). Therefore the matrix for $(\ ,\ )$ with respect to the orthogonal basis $\{e_0,e_1,\ldots, e_{d-1}\}$ has determinant  $(1-q)^{d(d-1)}$.  By (\ref{detW}) it follows that the matrix 
for $(\ ,\ )$ with respect to $\{ \omega, t(\omega) , t^2(\omega) , \ldots , t^{d-1}(\omega)\}$ has unit determinant. (Here we use that $\overline{1-q}=1-q^{-1}\sim 1-q$.) In other words, the form  $(\ ,\ )$ restricted to  $\mathcal {W}$ is unimodular. But then $\mathcal {W}$ must be equal to $\BS(S^1 \times S^1)$. This completes the proof.
\end{proof}

\begin{cor} \label{wc} If $H$ is a
 3-e-handlebody
 with boundary the
 e-surface
 of $\Sigma$ and $\Sigma$ has no colored points in the boundary, then $
\BS(\Sigma)$  is generated over $\mathcal{O}$ by elements represented by $\omega$-colored banded links  in $H.$ 
\end{cor}
\begin{proof} By the above theorem, each $e_i$ (in particular $e_1=z$) can be expressed as an $\mathcal{O}$-linear combinations of  the elements $t^j(\omega)$.
Therefore every mixed graph can be written as an $\mathcal{O}$-linear combination of $\omega$-colored banded links  in $H.$  The result now follows from Theorem~\ref{ml}.
\end{proof}

\section{Second integral basis in genus 1}\label{2nd}

Consider $K(d-1)$ in the notation of Section~\ref{Kauf}, now taking  $A= \zeta_{2p}.$ 
Let $\mathcal{V}$ denote its image in $\Vp(S^1 \times S^1).$ In other words  $\mathcal{V}$ is the  $\mathcal{O}$-submodule of $\Vp(S^1 \times S^1)$ 
generated by $\{ 1, v ,v^2 , \ldots , v^{d-1}\},$ where $v=(z+2)/(1+A)$. 
\begin{thm} \label {secondbasis} One has $\mathcal{V}=\BS(S^1 \times S^1).$ In particular, $\{ 1, v ,v^2 , \ldots , v^{d-1}\}$ is a basis for the  free module $\BS(S^1 \times S^1).$ 
\end{thm}

We refer to this basis as the $v$-basis of $\BS(S^1 \times S^1).$ We originally found it by the procedure outlined in Remark \ref{5.7}. Since $\tilde{\Gamma}(\Sigma)$ preserves $\BS(S^1 \times S^1)$, we have the following Corollary.

\begin{cor} $\mathcal{V}=\text{\em Span}_\BO\{1,v,v^2,\ldots,v^{d-1}\}$ is stable under the action of the mapping class group $\tilde{\Gamma}(S^1 \times S^1)$.
\end{cor}

\begin{rem} {\em  The mapping class group $\tilde{\Gamma}(S^1 \times S^1)$ is a central extension of $SL(2,\BZ)$. Its image in $GL(\Vp(S^1 \times S^1))$ is generated by $\kappa $ times the identity matrix (the central generator acts as multiplication by $\kappa$), the twist map $t$, and the so-called $S$-matrix. The entries of the $S$-matrix in the $e_i$-basis are well-known. One can therefore write down its entries in the $v$-basis (using the change of basis formulas in Lemma~\ref{changeofbasis}). The fact that these entries lie in $\BO$ is by no means obvious. We originally proved this fact using some  identities involving binomial coefficients.  The argument is similar to the proof that the $v$-basis is stable under the twist map $t$ given in Section~\ref{Kauf}, but considerably more complicated. We found proofs of these identities using Zeilberger's algorithm together with some identities from \cite{B} as above. In particular the {\tt Gosper} command in the
Mathematica package ``Fast Zeilberger'' (V 2.61) by Peter Paule and
Markus Schorn, \cite{Zeil} was used. As the proof we give below is much simpler, we omit the details of this computation.
}\end{rem}

\begin{proof}[ Proof of Theorem~\ref{secondbasis}]
\begin{lem} One has $\omega \in \BV~.$
\end{lem}
\begin{proof} Let $\lambda_i= -q^{i+1}- q^{-i-1}$.  
Recall \cite{BHMV1} that $e_i$ is an eigenvector with eigenvalue $\lambda_i$
for the endomorphism $c$ of $K (S^1 \times D^2,\BZ[A^{\pm}]) $ given by sending a skein in $S^1 \times D^2$ to the skein circled by a meridian.

 Let $\langle\ ,\ \rangle_H$ be the Hopf pairing ({\em i.e.} the symmetric bilinear form on $\Vp(S^1 \times S^1)$ which sends two elements $x,y$ to the bracket of the 
zero-framed
Hopf link with one component cabled by $x$, and the other component cabled by $y$). Then 
\begin{equation}\label{omip} \langle \omega,e_i\rangle_H = \begin{cases}
\langle\omega\rangle=\BD, &\text{if } i=0\\
\ \ \ 0  &\text{if } 1\leq i\leq d-1
\end{cases}\end{equation} 
Note that  $\langle z-\lambda_i,e_i\rangle_H=0$ for $1\leq i\leq d-1$, and $\langle z-\lambda_i,e_0\rangle_H=\lambda_0-\lambda_i$ (since $\langle z,e_0\rangle = \langle z\rangle =-q-q^{-1}=\lambda_0$). It follows that 
\begin{equation}\label{omi} \omega=\BD \prod_{i=1}^{d-1} \frac {z-\lambda_i} {\lambda_0-\lambda_i}\end{equation}
since the pairing $\langle\ ,\ \rangle_H$ is non-degenerate. (Note the similarity with the polynomials $Q_n$ of  \cite{BHMV1}.) 

Since $\BD\sim (1-q)^{d-1}$ and $\lambda_0-\lambda_i\sim (1-q)^2$ by Lemma~\ref{num}, it follows that 
\begin{equation}\label{omi2} \omega \sim \prod_{i=1}^{d-1} \frac {z-\lambda_i} {1-q}\sim \prod_{i=1}^{d-1} \frac {z-\lambda_i} {1+A}\end{equation} (where $\sim$ means equality up to multiplication by a unit). Now 
\begin{align*}
z-\lambda_i &= (z+2)-(2 +\lambda_i)\\
                        &= (z+2)-(1-q^{i+1})(1-q^{-i-1})\\
                        &= (z+2)+ u_i (1+A)^2
\end{align*}
where $u_i \in \BO.$  It follows that $$(z-\lambda_i)/(1+A) \in \text{Span}_\BO\{1,v\}~,$$  and so (\ref{omi2}) implies $\omega \in \BV~,$ proving the lemma.
\end{proof}

By Theorem~\ref{tK},  $K(n)$ hence $\BV$ is stable under the twist map $t$. It follows that  $$\mathcal{W}=\text{Span}_\BO\{\omega, t(\omega),\ldots, t^{d-1}(\omega)\} 
\subseteq \mathcal{V}~.$$  
Now recall from the proof of Theorem~\ref{firstbasis} that the matrix $W$ 
which expresses $\{\omega, t(\omega),\ldots, t^{d-1}(\omega)\}$ in terms 
of   $\{e_0,e_1,\ldots,e_{d-1}\}$ has determinant 
$\det W \sim
(1-q)^{-d(d-1)/2}$. 
Remembering $v=(z+2)/(1+A)$ and $1+A\sim 1-q$, it is easy to see that the same is true for the matrix which expresses $\{ 1,v,\ldots,v^{d-1}\}$ in terms of     $\{e_0,e_1,\ldots e_{d-1}\}$. Since $\mathcal{W}\subset\mathcal{V}$, it follows that actually $\mathcal{W}= \mathcal{V}$. By Theorem~\ref{firstbasis} we conclude  $\mathcal {V}=\BS(S^1 \times S^1)$. This completes the proof.
\end{proof}

By a $v$-colored banded link in a 3-manifold, we mean a banded link whose components are colored $v.$ As before this should be interpreted as the linear combination (superposition)  of the colored banded links that one obtains by expanding multilinearly. We note that $i$ parallel strands colored $v$ is the same as one strand colored $v^i.$

\begin{cor}\label{v}  If $H$ is a 
3-e-handlebody
 with boundary the 
e-surface
$\Sigma$ and $\Sigma$ has no colored points in the boundary, then $\BS(\Sigma)$  is generated over $\mathcal{O}$ by elements represented by $v$-colored banded links  in $H.$ 
\end{cor}

\begin{rem}{\em The matrix of the Hermitian form $(\ ,\ )_{S^1\times S^1}$ in the $v$-basis is easily computed. One has for $0\leq i,j\leq d-1$ 
\begin{align*}(v^i,v^j)_{S^1\times S^1}&=\BD \langle v^i,v^j\rangle_{S^1\times S^1} =  \langle v^{i+j},\omega \rangle_H
={(1+A)^{-(i+j)}} \langle (z+2)^{i+j},\omega \rangle_H\\
&=\frac {1} {(1+A)^{i+j}} \binom{2i+2j+2} {i+j} \frac 1 {i+j+1} \langle e_0,\omega \rangle_H\\
&=\frac {\BD} {(1+A)^{i+j}} \binom{2i+2j+2} {i+j} \frac 1 {i+j+1}
\end{align*}
Here we have used Lemma~\ref{changeofbasis} to express  $(z+2)^{i+j}$ in terms of the $e_n$, and then retained only the $e_0$ term. Indeed, the others are annihilated by the Hopf pairing with $\omega$ since $0\leq i+j\leq 2d-2$ (see (\ref{omip}) and remember that $e_{d+i}=e_{d-1-i}$ in $\Vp(S^1\times S^1)$).

It is instructive to 
check directly that the expression above 
lies in $\BO$ (use that $p$ divides the binomial coefficient 
$\binom{2i+2j+2} {i+j}$ if $d\leq i+j\leq 2d-2$).
}\end{rem}

\section{Integral basis in genus 2}\label{g2}

Let $\Si_2$ be a closed surface of genus $2$. In this section, we 
describe a basis for the module $\BS(\Si_2)$ and show that the Hermitian form $(\ ,\ )_{\Si_2}$ is unimodular.

Let $H_2$ be a regular neighborhood of the hand cuff graph $\Gtwo$ 
in $\mathbb{R}^3$. Then $H_2$ is a genus $2$ handlebody and by 
Corollary~\ref{v}, $\BS(\Si_2)$ is spanned by $v$-colored banded links 
in $H_2$. We think of  $H_2$ as
 $P_2 \times I$ where $P_2$ is a disk with two holes.

\begin{figure}[h]
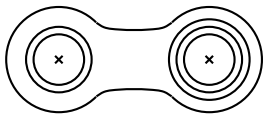

\begin{center}
\Ccurves
\caption{\label{fig0}The arrangement of curves $C_{2,3,1}$.} 
\end{center}
\end{figure}

The skein module  
$K(H_2,R)$ is free on the set of isotopy  classes  of collections  of nonintersecting essential simple closed curves in $P_2$.   
We refer to these isotopy classes as arrangements of curves. 
Such arrangements can be indexed by 3-tuples of 
nonnegative
integers. Let $C_{\alpha,\beta,\gamma}$ denote
the arrangement with  $\gamma$ parallel curves going around both holes, and within them $\alpha$ parallel curves  going around the left hole, and $\beta$ parallel curves  going around the right hole. See Figure \ref{fig0} for an 
example.

\begin{thm}\label{thg2} Let  $C_{\alpha,\beta,\gamma}(v)$  be the element of $\BS(\Si_2)$ obtained by coloring each curve of $C_{\alpha,\beta,\gamma}$ by $v=(z+2)/(1+A)$. Then the set 
   $$\{ C_{\alpha,\beta,\gamma}(v)\,|\,  0\le \gamma \le d-1, \quad
0\le \alpha,\beta  \le d-1-\gamma \} $$
 is a basis of $\BS(\Si_2).$  Moreover, the Hermitian form $(\ ,\ )_{\Si_2}$ is unimodular.
\end{thm}

Note that $C_{\alpha,\beta,\gamma}(v)$ lies in $\BS(\Si_2)$ because $v$ lies in $\BS(S^1\times S^1)$ by Theorem~\ref{secondbasis}.

\begin{proof} Let us first describe a basis of $\Vp(\Si_2)$ consisting of elements represented by colorings of the hand cuff graph \Gtwo . Let $G(i,j,k)$ be the element defined by the colored graph  
\begin{center}
\Gtwop
\end{center}
For this element to exist, $k$ must be even. Then the coloring 
is $p$-admissible if 
and only if
$\frac k 2 \leq i,j\leq p-2-\frac k 2$ (see \cite[Thm 1.15]{BHMV2}).
The standard basis of $\Vp(\Si_2)$ would be to take the $p$-admissible $G(i,j,k)$ with both $i$ and $j$ even. It is also possible to impose that one or both of $i,j$ be odd \cite[Thm 4.14]{BHMV2}. We will need a different basis where $i,j$ are allowed to be both  even and odd, but $\leq d-1$. This is given in the following Lemma.

\begin{lem}\label{lemG} The $G(i,j,k)$ with $k$ even in the range $[0,p-3]$, and both $i$ and $j$ in the range $[\frac k 2, d-1]$ (but not necessarily even), form a basis of $\Vp(\Si_2)$.
\end{lem}
\begin{rem}{\em 
Let $ \mathcal{G}$ be the basis  described in the above 
Lemma.
Let 
$ \mathcal{G}_k$
be the subset of elements of 
$ \mathcal{G}$
with middle arc colored $k$. The cardinality of 
$ \mathcal{G}_k$
is $(d-\frac k 2)^2.$ Thus we see directly that the cardinality of this basis
is  $  \sum_{j=0}^{d-1} (d-j)^2= \sum_{j=1}^{d} j^2= d (d+1)(2d+1)/6 .$
}\end{rem}

Lemma~\ref{lemG} could be proved using the methods of \cite{BHMV2}. Here we give a different, more direct proof.  

\begin{proof}[Proof of Lemma~\ref{lemG}] For $i$ in the range $[0,p-2]$ we let $i'=p-2-i$.  We claim that 
$$G(i,j,k)\sim G(i',j,k)\sim G(i,j',k)\sim G(i',j',k)$$ 
(where $\sim$ means equality up to multiplication by a unit in $\BO$). It is enough to prove that $G(i,j,k)\sim G(i',j,k).$ This is done in Figure \ref{fig1}. Note that if $i$ is even and $> d-1$ then $i'$ is odd and $\leq d-1$. Thus the basis of where all $i,j$ are even may be replaced by the basis of Lemma~\ref{lemG}.

\begin{figure}[h]
\begin{center}
\Ccalca \ \ \ \ $=$ \ \ \ \ \Ccalc \\
\hspace{2cm} $ = c_1$    \Ccalci \ \ \ $ = c_1 c_2$ \ \ \Ccalcai
\caption{\label{fig1} The proof that $G(i,j,k)\sim G(i', j, k)$. We have 
put $\tilde p= p-2$. In the first step, adding a loop colored $\tilde p$ 
doesn't change anything, since $e_{p-2}=e_0=1$. In the second 
(resp. third) step, we use the equation on the top (resp. bottom) 
of page 367 of \cite{MV}.
In the first equation, the sum reduces to only one term because 
 the other terms involve graphs which are not $p$-admissible. 
In the second equation, the coefficient is also {\em a priori}
a sum (coming from the expression for the tetrahedron coefficient
 in \cite[Theorem~2]{MV}).  But in fact here the sum
has only one term
(because  in the notation of \cite{MV} one has 
$\max(a_j)=p-2=\min(b_i)$ for this particular
tetrahedron coefficient).
The explicit formulas for $c_1$ and $c_2$ now 
show that both coefficients are products
of non-zero quantum integers 
or
their inverses
and therefore
are units
in $\BO$  (see Lemma~\ref{num}).}  
\end{center}
\end{figure}

\end{proof}

Let $ \mathcal{A}(v)$ be the set of the $v$-colored elements 
$C_{\alpha, \beta,\gamma}(v)$ claimed to be a basis in Theorem~\ref{thg2}, and let $ \mathcal{A}$ be the set of the uncolored
({\em i.e.} colored by $z=e_1$)
elements $C_{\alpha, \beta,\gamma}$ (in the same range for $\alpha, \beta,\gamma$).

\begin{lem} The set $ \mathcal{A}$ is a basis of $\Vp(\Sigma_2)$. 
Moreover, the 
basis
change from $ \mathcal{G}$ to $\mathcal{A}$ has determinant $\pm 1$.
\end{lem}
\begin{proof} Using the Wenzl recursion formula for the idempotents of the Temperley-Lieb algebra, one 
can expand the elements of $\mathcal{A}$ as $\mathcal{O}$-linear combinations of elements of
the graph basis
$\mathcal{G}.$ 
In fact, in the expansion of $C_{\alpha,\beta,\gamma}$, only those
$G(i,j,k)$ occur where $i\leq \alpha +\gamma$, $j\leq \beta +\gamma$,
and $k\leq 2\gamma$; moreover, $G(\alpha +\gamma,\beta +\gamma,
2\gamma)$ occurs with coefficient one. We can find orderings of
$\mathcal{A}$ and $\mathcal{G}$ so that the matrix which expresses
$\mathcal{A}$ in terms of $\mathcal{G}$ is triangular with ones on the
diagonal (use the lexicographical orderings where $\gamma$ resp. $k$
is counted first). This implies the Lemma.
\end{proof}

Let $ r=d (d+1)(2d+1)/{6}$ be the rank of  $\Vp(\Sigma_2)$. By Proposition~\ref{detb}, the matrix for $(\ ,\ )_{\Sigma_2} $ with respect to the orthogonal basis $\mathcal{G}$ has determinant $\sim (1-q)^{2 (d-1) r}$. The 
 preceding
 Lemma shows that the same holds true for the matrix for
 $(\ ,\ )_{\Sigma_2} $ with respect to $\mathcal{A}$.

Let $N$ denote the sum over $ \mathcal{A}$ of the number of curves appearing in each arrangement. The change of basis matrix for writing
$\mathcal{A}(v)$  in terms of $\mathcal{A}$ is again triangular and has determinant  $\sim (1-q)^{-N}.$ Thus the matrix  for $(\ ,\ )_{\Sigma_2} $ with respect to $\mathcal{A}(v)$ has determinant $\sim (1-q)^{2 (d-1) r -2N}$. The following Lemma~\ref{lemf} shows that this determinant is a unit. As in the genus one case (see the proof of Theorem~\ref{firstbasis}), we conclude that $\mathcal{A}(v)$ is a basis for $\BS(\Si_2)$ and that the form $(\ ,\ )_{\Sigma_2} $ is unimodular on $\BS(\Si_2)$.

\begin{lem} \label{lemf} $N=(d-1)r$. \end{lem}
\begin{proof}

To count $N, $ we write $\BA= \cup_{ 0\le \gamma \le d-1} \BA_\gamma,$ where 
\[ \BA_\gamma= \{ C_{\alpha,\beta,\gamma} |  0 \le \alpha,\beta  \le d-\gamma-1 \}~. \]

Note that $|\BA_\gamma|= (d-\gamma)^2.$ The total number of curves appearing in $\BA_\gamma$
is

 $$ \gamma (d-\gamma)^2 + \sum_{\alpha=0}^{d-\gamma-1}\sum_{\beta=0}^{d-\gamma-1} (\alpha +\beta) =  (d-1) (d-\gamma)^2~.$$
Thus each $\BA_\gamma$ contributes $d-1$ times its cardinality to the count. As 
$\sum_{\gamma =0}^{d-1} |\BA_\gamma| =r,$ we see that
$N=(d-1) r.$  
\end{proof}

This completes the proof of Theorem~\ref{thg2}.
\end{proof}

 \section{Non-Unimodularity 
}\label{non-uni}

Even without knowing an explicit basis of $\BS(\Si_g)$, it is possible to see that the form $(\ ,\ )_{\Si_g}$ is sometimes not unimodular. 

\begin{thm}\label{9.1} If $p \equiv 1 \pmod{4}$ and both the genus $g$ and the rank of $\Vp(\Si_g)$ are odd, then the form $(\ ,\ )_{\Si_g}$ is not unimodular on $\BS(\Si_g)$.
\end{thm}

For example, if $g=3$ and $p=5$ then the rank is $15$ and the form $(\ ,\ )_{\Si_3}$ is not unimodular on $\mathcal{S}_5(\Si_3)$.
\begin{rem} {\em
We used Mathematica \cite{Wo} to calculate 
the rank of $V_{p}(\Sigma_g)$
for small  $g$ using  the formula \cite[1.16(ii)]{BHMV2}.
We found that:
\begin{align*} \text{rank} \left( V_{4k+1}(\Sigma_3)\right) &=
 (1 /{45} )(3 \ k+32 \ k^2+120 \ k^3+200 \ k^4+192 \ k^5+128 \ k^6)
\\
\begin{split}
  \text{rank}   \left(V_{4k+1}(\Sigma_5)\right) &=
( 1/ {14175} )  (45 \ k + 864 \ k^2 + 6892 \ k^3 + 30184 \ k^4\\
&  +
      83760 \ k^5 + 172512 \ k^6 + 304896 \ k^7 + 458112 \ k^8\\ 
& + 542720 \ k^9 + 487424 \ k^{10} + 294912 \ k^{11} + 98304 \ k^{12})
\end{split} 
\end{align*}
 \noindent Thus the rank of $V_p(\Si_3)$ is odd if 
$p  \equiv 5 \pmod{8},$ and 
the form $(\ ,\ )_{\Si_3}$ is not unimodular 
in this case.
Similarly $(\ ,\ )_{\Si_5}$ is not unimodular if  
$p  \equiv 5 \pmod{8}. $  
}\end{rem}

        \begin{proof}[Proof of Theorem~\ref{9.1}]
The argument  relies on the following result of \cite{G}.
Assume $p \equiv 1 \pmod{4}$ and recall that $\BO=\BZ[\zeta_{4p}]$ in this case. Put $\BO^+=\BZ[\zeta_{p}]\subset \BO$. 
Let the Lagrangian assigned to $\Sigma_g$ be the kernel of the map on the 
first homology induced by the  inclusion of  $ \Sigma_g$ to $H_g$ and assign $H_g$  the weight zero. Then
$H_g$ is an { \em even}  (in the sense of  \cite{G}) morphism from $\emptyset$ to $\Sigma_g.$  Note that the quantum integers
$[n]$ for $1\le n \le p-1$ are units in $\BO^+.$

\begin{thm}\cite{G} If $p \equiv 1 \pmod{4}$ 
then $\BS(\Si_g)\simeq \BSplus(\Si_g)\otimes \BO$ where 
$\BSplus(\Si_g)\subset \BS(\Si_g)$ is a free $\BO^+$-module.
Moreover, one has $\mathcal {G}\subset\BSplus(\Si_g)$, where $\mathcal {G}$
is the colored graph basis of $\Vp(\Si_g)$ (see Proposition~\ref{detb}).
\end{thm}

The
matrix of $(\ ,\ )_{\Si_g}$ with respect to $\mathcal {G}$ has determinant $\BD^{gr}\sim (1-q)^{(d-1)gr}$ where $r$ denotes the rank of  $\Vp(\Si_g)$.  Let $\mathcal {B}$ be a basis of the free  $\BO^+$-module $\BSplus(\Si_g)$, and let $D$ be the determinant of the matrix expressing $\mathcal {B}$ in terms of $\mathcal {G}$. The matrix of $(\ ,\ )_{\Si_g}$ with respect to the basis $\mathcal {B}$ has determinant $\sim \Delta$, where 
\begin{equation}\label{Delta} 
\Delta =   D \overline D 
(1-q)^{(d-1)gr}~.
\end{equation}  If the form  is unimodular, $\Delta$ must be a unit
in $\BO$, 
and since $\Delta$ lies in $\BO^+$, it must be a unit in $\BO^+$. But $1-q$ is a self-conjugate prime in 
$\BO^+=\BZ[q]=\BZ[\zeta_{p}]$, and since 
$D^{-1}$
lies in $\BO^+$ as well, 
$\Delta$ can be a unit only
if $(d-1)gr$ is even. Thus one of $g$ and $r$ must be even (since $d-1=(p-3)/2$ is odd in our situation). This completes the proof.
\end{proof}
\begin{rem}\label{9.4} {\em The use of the $\BO^+$-module $\BSplus(\Si_g)$ can in general not be avoided in this argument. Here is why. Recall that $1-q$ splits in $\BO=\BZ[\zeta_{4p}]$ as 
the product of two conjugate prime ideals $\mathfrak{p}$ and $\overline{\mathfrak{p}}$. If $\mathfrak p$ is principal (this happens for example if $p=5$), then there exists $D\in \BO$ such that 
the number $\Delta$ defined as in (\ref{Delta}) is a unit
even when $(d-1)gr$ is odd. 
Of course, such a $D$ does not exist in $\BO^+$.}\end{rem}

\begin{rem} {\em  If we assign extra structure to $\Si_g$ and $H_g$ as described above in the proof of 9.1, then   $\BO^+$ linear combinations of banded links in $H$ represent  elements in $\BSplus(\Si).$
 Moreover the bases described in Sections~\ref{1st},  \ref{2nd}, \ref{g2} for $\BS(S^1 \times S^1), $ and  $\BS(\Si_2)$  are actually bases for  
$\BSplus(S^1 \times S^1), $  and $\BSplus(\Si_2).$  There are also plus versions of Theorem
\ref{ml} and Corollaries \ref{wc} and \ref{v}.}\end{rem}

\begin{rem} {\em
When restricted to $\BSplus(\Sigma_g)$, the Hermitian
form $(\ ,\ )_{\Sigma_g}$ does not take values in $\BO^+,$ if $g$  is odd.
This follows from
the proof of
Proposition~\ref{detb}, since $\BD \not\in \BO^+$.
In the next section,
we will use the 
sesquilinear form
\[ (\ ,\ )^+_{\Sigma_g}: \BSplus(\Sigma_g)
\times
\BSplus(\Sigma_g) \rightarrow \mathcal{O}^+\]
obtained by multiplying the form $(\ ,\ )_{\Sigma_g}$
by $i^{\varepsilon (g)},$ where ${\varepsilon (g)}$ is zero or one 
accordingly   
as $g$ is even, or odd.
This form takes values in $\BO^+$ since $i \BD \in \BO^+$.
}\end{rem}

\section{Genus three at the prime five}\label{g3}

In genus $g\geq 3$, 
one can also try
to find a set  of 
banded links
in a  handlebody so that one obtains a basis of $\BS(\Si_g)$
 by cabling each curve
component
 with $v=(z+2)/(1+A)$. 
This is suggested by 
 Corollary \ref{v} and the fact that $\BS(\Si_g)$ is a 
free $\BO$-module \cite{G}. 
In fact, we now find such a set of 
links 
giving a basis for 
$\mathcal{S}^+_5(\Si_3)$ (and therefore also for $\mathcal{S}_5(\Si_3)$)
by adapting the above procedures. These links are described by 
arrangements of curves in
a 
 thrice punctured disk.  Although the
 Hermitian form and the related 
$\BO^+$-valued 
sesquilinear form are not unimodular, 
in this particular situation
they are nearly so, and this is essential for our argument.
It seems more difficult to find an explicit collection of banded links with this property
for $\mathcal{S}_p(\Si_3)$ for $p>5,$ and for  $\mathcal{S}_p(\Si_g)$ for $g>4.$
We plan to return to this question elsewhere.

We think of  
the handlebody
$H_3$ as
 $P_3 \times I$ where $P_3$ is a disk with three holes. We give $H_3$ weight zero. We equip $\Si_3$ with the Lagrangian given by the kernel of the map induced on the first homology by the inclusion of  $\Si_3$ in $P_3 \times I.$ 

\begin{figure}[h]
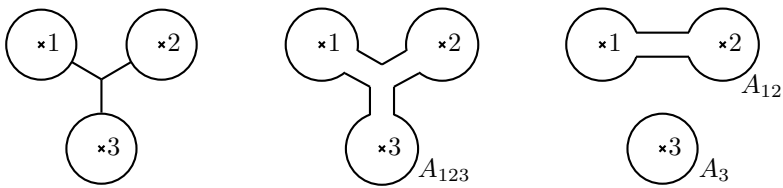

\begin{center}
\Gtroisi \ \ \ \ \ \ \ \ \Gtroisii \ \ \ \ \ \ \ \ \Gtrois 
\caption{\label{fig3}The graph $G$ and the arrangements of curves $A_{123}$ and $A_{12}A_3$.} 
\end{center}
\end{figure}

Consider the set of 15 arrangements of curves in $P_3$ 
\begin{align*}\BA=\{A_\emptyset,&A_1,A_2,A_3, A_1A_2,A_2A_3,A_3A_1, 
A_1A_2A_3,\\
 &A_{12},A_{23},A_{13}, A_{12}A_3, A_{23}A_1, A_{31}A_2, A_{123}\}~.
\end{align*}
Here, $A_\emptyset$ is the empty arrangement, $A_i$ (resp. $A_{ij}$, resp. $A_{123}$) is a 
curve of  the shape  pictured in Figure~\ref{fig3} around
just
 the $i$-th hole (resp.
around both 
the $i$-th and $j$-th hole, resp. 
around all three holes), and the multiplicative notation $A_\alpha A_\beta$ means disjoint union of $A_\alpha$ and $A_\beta$. See Figure~\ref{fig3} for two examples. Note that the total number of curves in $\BA$ is $22$.

\begin{thm}\label{thg3} The set 
$\BA(v)=\{A_\emptyset(v)=A_\emptyset,A_1(v),A_2(v),\ldots\}$ 
consisting of the curve arrangements in $\BA$ colored $v$
 is a basis of $\mathcal{S}^+_5(\Si_3), $  and thus also a basis for $\mathcal{S}_5(\Si_3).$ 
\end{thm}
Note that it follows in particular that the $\BO^+$-span of $\BA(v)$ is stable under the action of the 
index two subgroup of even morphisms in the
mapping class group $\tilde{\Gamma}(\Si_2)$.

\begin{lem} The set $\BA$ (where its elements are considered as planar banded links in $H_3$) is a basis of $V_5(\Si_3)$.
\end{lem}
\begin{proof}
By the proof of Lemma \ref{lemG}, we can find a graph basis 
$\mathcal{G}$
for $V_5(\Si_3)$ by 5-admissible colorings of the graph $G$ in Figure~\ref{fig3},
where the loops are colored zero or one and the non-loop edges are colored zero or two.  
If a non-loop edge is colored two, then the loop at the end of the edge 
must be colored one. Also, the number of non-loop edges colored two must be
zero, two, or three. 
This summarises  5-admissibility in this case. There are 15 such colorings. 
Again using Wenzl's recursion formula, there is a triangular change of 
basis with ones on the diagonal 
from the basis $\mathcal{G}$ to the set 
$\mathcal{A}$
which is therefore also a basis of  $V_5(\Si_3)$.
\end{proof}

\begin{proof}[Proof of Theorem~\ref{thg3}]
Recall that $\BA(v)$ consists of the $15$ elements of $\BSfplus(\Sigma_3)$ 
obtained by replacing each of the 22 
curves in $\BA$ with $v=(z+2)/(1+A)$.  Again there is a triangular change 
of basis 
matrix from $\BA$ to $\BA(v)$.  Therefore the elements of $\BA(v)$
span $V_5(\Si_3)$ and hence are linearly independent over $\BO^+$.
Consider the inclusion 
\begin{equation}\label{strict}
\text{Span}_{\BO^+}\BA(v) \subset \BSfplus(\Sigma_3)~.
\end{equation} The matrix 
for $(\ ,\ )^+_{\Sigma}$ with respect to  $\mathcal{A }(v)$ has 
determinant 
$\sim (1-q)^{3\cdot 15-2 \cdot 22}=1-q.$ 
Since 
$1-q$ is a prime in  $\BO^+ $  and $\BSfplus(\Sigma_3)$ is also 
a free $\BO^+$-module,
we conclude that the inclusion (\ref{strict}) 
cannot be strict. Thus $\mathcal{A }(v)$ is a basis for $\BSfplus(\Sigma_3).$
\end{proof}

\begin{rem}{\em Theorem \ref{thg3} remains true if we replace $v$ by $\omega$ throughout.
The same proof works.}
\end{rem}

\begin{rem} {\em As in Remark~\ref{9.4}, it is crucial for this argument 
to use $\BSfplus(\Si_3)$ rather than $\BSf(\Si_3)$, since for $p=5$ there exists $a\in \BO=\BZ[\zeta_{20}]$ such that $1-q=a\overline{a}$.}\end{rem}

\begin{rem}{\em Kerler has announced in \cite{Ke} a construction of integral
bases for the Reshetikhin-Turaev $SO(3)$
TQFT at the prime $p=5$ for any genus.}\end{rem}

\section{A divisibility result for the Kauffman bracket} 

In this final section, we let $A$ again be an indeterminant. The fact 
that $v=(z+2)/(1+
\zeta_{2p}
)$ lies in $\BS(S^1\times S^1)$ for all odd primes $p$ 
has the following application to the Kauffman bracket $\langle\ \rangle$ 
of banded links in $S^3$.

\begin{thm} \label{th11} Let $L$ be a banded link in $S^3$ with $\mu$ components.
 Let $L(z+2)$ denote this link colored $z+2.$
Then the  Kauffman Bracket 
$\langle L(z+2) \rangle \in \BZ[A^{\pm}]$ is divisible by 
$(1+A)^{\mu}$. 
\end{thm}
Here the Kauffman bracket is 
normalized
so that the bracket of the empty 
link is $\langle \emptyset\rangle=1$. Note that 
$$\langle L(z+2)\rangle = \sum_{L'\subset L} 2^{\mu-\mu(L')} 
\langle L'\rangle~,$$ where the sum is over all sublinks 
$L'$ of $L$, and $\mu(L')$ denotes the  number of components of $L'$.

\begin{proof}[Proof of Theorem \ref{th11}] When we evaluate the Kauffman 
bracket $\langle J \rangle$ of a banded link $J$ in $S^3$ at $A=\zeta_{2p}$, 
we obtain the quantum invariant of the pair $(S^3,J)$ (where $S^3$ is given the weight zero) in the 
normalization
$$\langle J \rangle\vert_{A=\zeta_{2p}}=I_p(S^3,J)=\frac {\langle (S^3,J) \rangle} {\langle S^3 \rangle}=\BD {\langle (S^3,J) \rangle} ~.$$
This
normalization
$I_p(M,J)$ of the quantum invariant is precisely the 
one 
which is always an algebraic integer \cite{Mu1,MR} and which is at the basis of the integral cobordism functors $\BS$.

Let $f(A)$ denote the Kauffman bracket  $\langle L(z+2)\rangle \in \BZ[A^\pm].$ 
Since 
 $v= (z+2)/(1+
\zeta_{2p}
)\in S_p(S^1 \times S^1),$
we have $I_p(S^3,L(v)) \in \BZ[\zeta_{2p}],$ for every odd prime $p.$
Thus 
\begin{equation}
 \label{divisi}
f(\zeta_{2p}) = I_p(S^3,L(z+2)) \in (1+\zeta_{2p})^\mu \BZ[\zeta_{2p}],
\end{equation} for every odd prime $p.$ 

Now recall the following elementary Lemma (see \cite[Lemma 5.5]{Mu1} and note
 that  $-\zeta_{2p}$ is a primitive $p$-th root).

\begin{lem} 
Suppose $f(A) \in {\BZ} [A^\pm].$ Let $f^{(k)}(A)$ denote the $k$-th derivative of $f(A)$. Assume $0\leq \mu<p$ where $p$ is prime. Then 
 $f(\zeta_{2p})\in {\BZ} [\zeta_{p} ]$ is divisible by $(1+\zeta_{2p})^\mu$ if and only if $f^{(k)}(-1)\equiv 0\pmod{p}$ for every $0\le k< \mu.$
\end{lem}

By this lemma,  (\ref{divisi}) implies $f^{(k)}(-1)\equiv 0\pmod{p}$ for each $0\le k< \mu,$ provided $p$ is larger than $\mu.$ Since there are infinitely many such primes, it follows that $f^{(k)}(-1)=0$ for each $0\le k< \mu.$ 
But this means that $ (1+A)^{\mu}$ divides $f(A).$
  \end{proof}

\begin{cor} \label{Kcor} If $L$ is as in the theorem, then the  Kauffman Bracket 
$\langle L(z+[2]) \rangle \in \BZ[A^{\pm}]$ is also divisible by 
$(1+A)^{\mu}$.
\end{cor}
\begin{proof} This follows immediately from the fact that 
$2-[2]=2-A^2-A^{-2}=(1-A^2)(1-A^{-2})$ is divisible by $1+A.$
\end{proof} 

\begin{rem}{\em In a similar way, 
Theorem \ref{tK} remains true if we replace
$v$ with $\hat v= (z+ [2])/(1+A)$
in the definition  of $K(n).$  Similarly  Theorems \ref{secondbasis}, \ref{thg2}, \ref{thg3} remain true if we replace
$v$ by $\hat v.$
}\end{rem}

\begin{rem}{\em  Theorem \ref{th11} can also be proved by computing
the Kauffman bracket from the (framed)
Kontsevich integral via an appropriate weight system.   Actually this proof is
an
adaptation of an argument going back
to Kricker and Spence  \cite[Proof of Thm. 2]{KS}, but they only 
considered algebraically split links. Previously Ohtsuki 
\cite[Prop. 3.4]{OhCambridge} had obtained a stronger  
divisiblity result for $\langle L(z+[2]) \rangle$ using quantum groups, for algebraically split links satisfying some extra conditions. Later, Cochran and 
Melvin  generalised the Kontsevich integral argument, and their result \cite[Theorem 2.5]{CM}
contains Corollary \ref{Kcor} for zero-framed links. (The results of \cite{OhCambridge,KS,CM} are 
stated in terms of the Jones polynomial,  but it is well-known that the Jones polynomial and the Kauffman
bracket are equivalent.) However, the restriction
to zero framing is not really necessary 
(although a small additional argument
is needed).
We will not give details of
this alternative proof here, as the techniques are completely different
from the ones in the present paper. 
}\end{rem}

\end{document}